\documentclass[12pt]{article}

\usepackage[latin1]{inputenc}
\usepackage{amsfonts}
\usepackage{amsmath}
\usepackage{amsfonts}
\usepackage{latexsym}
\usepackage{amssymb}
\usepackage{dsfont}

\newtheorem{fed}{Definition}
\newtheorem{teo}[fed]{Theorem}
\newtheorem{lem}[fed]{Lemma}
\newtheorem{cor}[fed]{Corollary}
\newtheorem{pro}[fed]{Proposition}
\newtheorem{rem}[fed]{Remark}

\newtheorem{exa}[fed]{Example}
\newtheorem{exas}[fed]{Examples}

\def\proof{\noindent {\it Proof. }}
\def\In{\mathbb {I} _n}
\def\IM{\mathbb {I} _m}
\def\mJ{{\mathbb J}}
\def\mS{{\mathbb S}}

\def\suml{\sum\limits}
\def\QEDP{\tag*{\QED}}

\def\ov{\overline}
\def\bce{\begin{center}}
\def\ece{\end{center}}
\def\ds{\displaystyle}
\def\noi{\noindent}
\def\cF{\mathcal F}
\def\QED{\hfill $\square$}
\def\EOE{\hfill $\triangle$}
\newcommand{\peso}[1]{ \quad \text{ #1 } \quad }
\def\uno{\mathds{1}}
\def\bm{\left(\begin{array}}
\def\em{\end{array}\right)}
\def\ben{\begin{enumerate}}
\def\een{\end{enumerate}}
\def\barr{\begin{array}}
\def\earr{\end{array}}
\def\igdef{\ \stackrel{\mbox{\tiny{def}}}{=}\ }
\def\k{\mathbf{k}}
\def\w{\mathbf{w}}
\def\v{\mathbf{v}}
\def\eps{\varepsilon}
\def\N{\mathbb{N}}
\def\R{\mathbb{R}}
\def\C{\mathbb{C}}
\def\cI{\mathcal{I}}
\def\cH{\mathcal{H}}
\def\cK{\mathcal{K}}
\def\cP{\mathcal{P}}
\def\cR{{\cal R}}
\def\cS{{\cal S}}
\def\cM{{\cal M}}
\def\cN{{\cal N}}
\def\cV{{\cal V}}
\def\cU{{\cal U}}
\def\cW{{\cal W}}
\def\ese{\mathcal{S}}
\def\ene{\mathcal{N}}
\def\vacio{\varnothing}
\def\orto{^\perp}
\def\inc{\subseteq}
\def\sii{ if and only if }
\def\inv{^{-1}}
\def\rai{^{1/2}}
\def\api{\langle}
\def\cpi{\rangle}
\def\beq{\begin{equation}}
\def\eeq{\end{equation}}
\def\pausa{\medskip\noi}
\def\SP{L(m,\k,d)}
\def\RS{\cR \cS (m,\k,d)}
\def\RSI{\cR \cS _\cI (m,\k,d)}
\def\PRO{\cP (m,\k,d)}
\def\RSV{\cV= \{V_i\}_{i\in \, \IN{m}}}
\def\RSW{\cW= \{W_i\}_{i\in \, \IN{m}}}
\def\RSS{\cS= \{S_i\}_{i\in \, \IN{m}}}

\def\H{{\cal H}}

\def\uni{{\cal U}(\H )}
\def\glh{Gl(\cH)}
\def\coma{\, , \, }
\def\fori{for every $i \in \IM\,$}
\def\subim{_{i\in \,\IN{m}}\,}
\def\FS{\cN_w  = (w_i\, ,\, \cN_i)_{i\in \IM}}
\def\la{\lambda}

\DeclareMathOperator{\Preal}{Re}

\DeclareMathOperator{\tr}{tr}
\DeclareMathOperator{\gen}{span}

\newcommand{\pint}[1]{\displaystyle \left \langle #1 \right\rangle}
\newcommand{\hil}{\mathcal{H}}
\newcommand{\op}{L(\mathcal{H})}
\newcommand{\opk}{L(\mathcal{K})}
\newcommand{\lhk}{L(\mathcal{H} , \mathcal{K})}
\newcommand{\lkh}{L(\mathcal{K} , \mathcal{H})}
\newcommand{\posop}{L(\mathcal{H})^+}
\newcommand{\cene}{\mathbb{C}^n}
\newcommand{\mat}{\mathcal{M}_n(\mathbb{C})}

\newcommand{\matu}{\mathcal{U}(n)}
\newcommand{\matpos}{\mat^+}
\newcommand{\matinv}{\mathcal{G}\textit{l}\,(n)}
\newcommand{\matinvd}{\mathcal{G}\textit{l}\,(d)}
\newcommand{\matrec}[1]{\mathcal{M}_{#1} (\mathbb{C})}
\newcommand{\IN}[1]{\mathbb {I} _{#1}}

\begin{document}

\title{Robust dual reconstruction systems  and fusion frames}


\author{Pedro G. Massey, Mariano A. Ruiz  and Demetrio Stojanoff\thanks{Partially supported by CONICET 
(PIP 5272/05) and  Universidad de La PLata (UNLP 11 X472).} }
\author{P. G. Massey, M. A. Ruiz and D. Stojanoff \\ {\small Depto. de Matem\'atica, FCE-UNLP,  La Plata, Argentina
and IAM-CONICET  
}}
\date{}
\maketitle

\begin{abstract}
We study the duality of reconstruction systems, which are $g$-frames in a finite dimensional setting. These systems allow redundant linear encoding-decoding schemes implemented by the so-called dual reconstruction systems.
 We are particularly interested in the projective reconstruction systems that are the analogue of fusion frames in this context. Thus, we focus on dual systems of a fixed projective system that are optimal with respect to erasures of the reconstruction system coefficients involved in the decoding process. We consider two different measures of the reconstruction error in a blind reconstruction algorithm.
 We also study the projective reconstruction system that best approximate an arbitrary reconstruction system, based on some well known results in matrix theory. Finally, we present a family of examples in which the problem of existence of a dual projective system of a reconstruction system of this type is considered.
\end{abstract}

\pausa
{\bf keywords:} Reconstruction systems, dual  reconstruction systems, erasures, fusion frames.

\pausa
{\bf subclass:} 42C15 and 15A60.

\tableofcontents

\section{Introduction}

Many researchers have recently studied the problem of designing 
finite frames 
for the reconstruction of a signal $x$ in the absence of a number of (missing or delayed) coefficients (see \cite{Bod,Pau,BodPau,[CasKu2],HolPau,GK,P,Stro}). In this context, the goal is to construct finite frames which minimize the (blind) reconstruction error of a signal  when a number $k$ of frame coefficients are ``erased''.

In \cite{HanLo}, the authors study the following problem: given a {\it fixed} frame $\cF=\{f_i\}_{i=1}^m$ for $\C^d$, find the alternate dual frames that minimize a measure of the (blind) reconstruction error when $k$ frame coefficients are erased.
 Thus, unlike the previous setting,  the interest is centered in a subset of the set of alternate dual frames of a fixed frame. As they notice, the canonical dual is not always the best choice despite the fact that it is associated to the Moore-Penrose pseudoinverse. 
Nevertheless, one of the main results in \cite{HanLo} exhibits conditions on the fixed frame $\cF$ which assure that the only optimal dual of $\cF$ - in the sense previously described -  is its canonical dual.

Our research is in the vein of \cite{HanLo} but our interest is to set it in the context of fusion frames. However, since the definition of fusion frames depends strongly on the subspaces of the frame, the concept of dual fusion frame is not easy to handle. Therefore, we identify fusion frames with a subset of a broader class, namely the Reconstruction Systems (RS) (finite dimensional $g$-frames, using the terminology of W. Sun \cite{GFsun}). Hence we fix a {\it projective} RS $\cV$  (i.e. the analogue of fusion frames in the context of RS's), and we search for conditions that allow us to describe the dual RS's of $\cV$ that minimize some measure of the (blind) reconstruction error.

 As in  \cite{HanLo}, our results show the existence of a unique dual RS (which is not necessarily projective) that is optimal for the erasure of 1 packet of coefficients. As in previous works in the subject, the error is measured in terms of the Frobenius norms of the so-called ``packet lost operators'' (the  equivalent notion of the error operators in \cite{Bod,HolPau,HanLo}). We present two different ways to perform this measure: the worst case error (WCE), i.e. the largest packet lost operator norm among all possible erasures, and the 2-error which is the euclidean norm of the vector of all packet lost operator norms. Each of these errors bound the norm of the error when a vector $x$ is reconstructed with a blind reconstruction strategy, assuming that a frame coefficient is missing.

 For the WCE, following the approach of Han and Lopez, we impose some conditions on the norms involving the canonical dual elements to assure that the canonical dual is the only optimal alternate dual for one erasure. For the 2-error we show that, for projective RS's, there is always a unique optimal dual.  In both cases, our results include a detailed description of the optimal dual RS.

The fact that the optimal alternate duals of the fixed projective RS $\cV$ considered above need not be projective
motivates the study of best projective approximations of an arbitrary RS. We consider this problem in terms of well known results in matrix theory.

On the other hand, assuming that the erased packets of coefficients correspond to a known subset of indexes of the RS $\cV$,  we find conditions which assure that the RS corresponding to the remaining set of indexes is also a RS. Notice that this last fact allows to have perfect reconstruction with suitable dual RS's.
We have included a section of examples in which we also consider the problem of the existence of projective duals of a given projective RS.

In order to put our results in perspective and to keep the text self contained we have reviewed some known results concerning the structure of optimal duals  in the context of fusion frames and $g$-frames, including short proofs in some cases. 
It is worth pointing out that some results in $g$-frames theory can be obtained from their analogues in frame theory, using the fact that $g$-frames can be considered as groupings of vectors of a frame. We have used this kind of argument whenever it was possible. Nevertheless, when an extra structure is imposed to the RS's new difficulties may arise; such is the case of projective RS's (i.e. our analogue of fusion frames), where there are orthogonal relations between subsets of vectors, or injective RS's.  Indeed, as far as we know Theorem \ref{unico}  is new even in the vector frame setting. In addition, by the nature of its formulation, the approximation by projective RS's and the existence of projective dual RS's of a fixed projective RS are intrinsic problems to the RS's setting.

The paper is organized as follows: In Section \ref{basic} we recall
some basic facts about general reconstruction systems and
we fix some of the terminology used throughout the paper.
In Section \ref{sec erasures} we study the optimal dual systems
for erasures of vector coefficients.
In Section \ref{rem erasu lb} we consider the problem of perfect reconstruction when packets corresponding to a fixed set of indices of a RS $\cV$ are lost.
In Section \ref{Approx} we describe the  projective RS which is nearest to a fixed RS.
In Section \ref{Opti} we give several examples of group RS's, and we consider the problems of existence of projective dual RS's of a given projective RS.

\subsection*{Notations.}
Given $\cH \cong \C^d$  and $\cK \cong \C^n$, we denote by $\lhk $
the space of linear operators $T : \cH \to \cK$.
Given an operator $T \in \lhk$, $R(T) \inc \cK$ denotes the
image of $T$, $\ker T\inc \cH$ the null space of $T$ and $T^*\in \lkh$
the adjoint of $T$. If $d\le n$ we say that $U\in \lhk$ is an isometry
if $U^*U = I_\cH\,$. In this case, $U^*$ is called a coisometry.
We denote by $\cI(d\coma n)$ the set of all isometries in $\lhk$.

\pausa
If $\cK = \cH$ we denote by $\op = L(\cH \coma \cH)$,
by $\glh$ the group of all invertible operators in $\op$,
 by $\posop $ the cone of positive operators and by
$\glh^+ = \glh \cap \posop$.

\pausa
If $T\in \op$, we  denote by
$\sigma (T)$ the spectrum of $T$, by rk $T $  the rank of $T$,
and by $\tr T$ the trace of $T$.
Given $m \in \N$ we denote by $\IM = \{1, \dots , m\} \inc \N$ and
$\uno = \uno_m  \in \R^m$ denotes the vector with all its entries equal to $1$.

\pausa
On the other hand, $\matrec{n,m}$ denotes the space of complex $n\times m$ matrices.
If $n=m$ we write $\mat = \matrec{n,n}$,
$\matinv$ the group of all invertible elements of $\mat$, $\matu$ the group
of unitary matrices,
 $\matpos$ the set of positive semidefinite
matrices, and $\matinv^+ = \matpos \cap \matinv$.

\pausa
If $W\inc \cH$ is a subspace we denote by $P_W \in \posop$ the orthogonal
projection onto $W$, i.e. $R(P_W) = W$ and $\ker \, P_W = W^\perp$.
For vectors on $\cene$ we shall use the euclidean norm, but for matrices $T\in \mat$, we shall use both
\begin{enumerate}
\item The spectral norm $\|T\| = \|T\|_{sp}= \max\limits_{\|x\|=1}\|Tx\|$.
\item\label{item 2}   The Frobenius norm $\|T\|_{_2} = (\tr \, T^*T )\rai =
\big( \, \suml_{i,j \in \In } |T_{ij}|^2 \, \big)\rai$. This norm  is induced by the inner product
$\api A,\ B\cpi=  \tr \, B^*A \,$,  for $A, B \in \mat$.

\end{enumerate}
\pausa

\section{Basic framework of reconstruction systems}\label{basic}

In what follows we consider $(m,\k,d)$-reconstruction systems (see for example \cite{P,MRS2,GFsun}), which
are more general linear systems than those considered in
\cite{BF,Bod,Pau,BodPau,HolPau} and
\cite{MR}, that also have an associated reconstruction algorithm.
\begin{fed}\label{defi recons}\rm
Let $m, d \in \N$ and $\k = (k_1 \coma \dots \coma k_m)   \in \N^m $.
\ben
\item We denote by $
\cK = \cK_{m\coma \k} \igdef  \bigoplus_{i\in \, \IN{m}}  \ \C^{k_i} \ .
$
Sometimes we shall write each direct summand by $\cK_i = \C^{k_i}\,$.

\item Given a space $\cH \cong \C^d$ we denote by
$\SP \igdef \bigoplus_{i\in \IM} L(\cH \coma \cK_i) \cong \lhk$.
A typical element
of $\SP$ is a system $\RSV $ such that
each $V_i \in L(\cH \coma \mathcal K_i)$.
\item
A family $\cV= \{V_i\}_{i\in \, \IN{m}}\in \SP$  is an
$(m,\k,d)$-reconstruction system (RS)
for $\cH$ if
$$
\barr{rl}
S_\cV & \igdef \sum_{i\in \, \IN{m}} V_i^*V_i  \in \glh^+\, ,
\earr$$
i.e., if $S_\cV$ is  invertible and positive. $S_\cV$ is called the {\bf RS operator} of $\cV$. In this case,
the $m$-tuple  $\k = (k_1 \coma \dots \coma k_m)   \in \N^m $ satisfies that
$\tr \k \igdef \sum_{i\in \, \IN{m}} k_i \ge d$. 
\item
The bounds of $\cV$ as a RS are the numbers 
$A_\cV = \la_{\min} (S_\cV)$ and $B_\cV = \|S_\cV\|_{sp}\,$. Observe that
$$
\barr{rl}
A_\cV\,\|x\| ^2 &\le \api S_\cV\, x\coma x\cpi =  
\sum_{i\in \, \IN{m}} \|V_i \, x\|^2 \le B_\cV\,\|x\| ^2 
\peso{for every} x\in \cH\ . \earr
$$
\item
We shall denote by $\RS$ the set of all $(m,\k,d)$-RS's for $\cH\cong \C^d$.

\item
The system $\RSV \in \RS$  is said to be {\bf injective} if
$V_i^* \in L(\cK_i\coma \cH)$ is injective (equivalently, if $V_i\,V_i^* \in 	Gl(\cK_i)\,$) \fori .

We shall denote by  $\RSI$ the set of all injective elements of $\RS$.
\item
The system $\cV $  is said to be {\bf projective} if
there exists a sequence  $\v= (v_i)_{i \in \IN{m}}\in \R_+^m $
of positive numbers, the
weights of $\cV$, such that
$$
V_i \, V_i^* = v_i^2 \, I_{\cK_i} \ , \peso{for every} i \in \IN{m} \ .
$$
In this case, the following properties hold:
\ben
\item The weights can be computed directly, since each  $v_i = \|V_i\|_{sp}\ $.
\item Each $V_i = v_i U_i$ for a  coisometry $U_i \in L(\cH \coma \cK_i)$.
Thus $V_i^*V_i=v_i^2\,
P_{R(V_i^*)} \in \posop$ \fori.

\item
Observe that in this case
$S_\cV = \sum_{i\in \, \IN{m}} v_i^2 \,P_{R(V_i^*)}$ as in fusion frame theory.
\een
We shall denote by  $\PRO$ the set of all projective elements of $\RS$.

\item
The {\bf analysis} operator of the system $\cV$ is defined by
$$
\barr{rl}
T_\cV \ : \ \mathcal H\rightarrow
\cK  & =
 \bigoplus_{i\in \, \IN{m}}  \cK_i\ \ \text{given by } \ \ T_\cV\, x
 = (V_1 \, x \coma \dots \coma V_m\, x)
\ , \peso{for} x\in \cH \ .\earr
$$
\item
Its adjoint  $T_\cV^*$ is called the {\bf synthesis} operator of the system $\cV$, 
and it satisfies that
$$
\barr{rl}
T_\cV^* : \cK = \bigoplus_{i\in \, \IN{m}}  \cK_i \rightarrow\mathcal H  
&\peso{is given by}
T_\cV ^* \, \big(\,  (y_i)_{i\in \, \IN{m}}\big) =\sum_{i\in \, \IN{m}} V_i^*\, y_i \ .
\earr
$$
Using the previous notations and definitions we have that $S_\cV = T_\cV^*\  T_\cV\, $. \EOE
\een
\end{fed}

\begin{exas}[Vector and fusion frames as RS's]\label{ejems} \

\pausa \rm
1.
As it was mentioned earlier,  RS's arise from usual vector frames
by grouping together the elements of the frame.
Therefore, it is natural to expect that in the case
$\k = \uno_m\,$, the set $\RS$ can be identified with the set of $m$-vector frames
for $\cH \cong \C^d$.

\pausa
Indeed, let $\cF=\{f_i\}_{i\in \, \IN{m}}\in \H^m $.
For $i\in \IN{m}$ consider $V_i:\H\rightarrow \C$ given by $V_i(x)=\langle x, f_i\rangle$ for every $x\in \cH$.
Let $\cV_\cF=\{V_i\}_{i\in \IN{m}}$ and notice that
$$
S_{\cV_\cF}=\sum_{i\in \IN{m}} V_i^*V_i =\sum_{i\in \IN{m}} \langle \,\cdot\,,f_i\rangle \,f_i=S_\cF \ .
$$
Thus $\cF$ is a frame for $\cH$  if and only if $\cV_\cF\in \mathcal{RS} (m,\uno,d)$.
Actually,
$\mathcal{RS} (m,\uno,d) =\mathcal{P} (m,\uno,d)$ because
every functional is a multiple of a coisometry.
Moreover,  $T_{\cV_\cF}:\H\rightarrow \oplus_{i\in \, \IN{m}} \C=\C^m$ is the usual
analysis operator of $\cF$. On the other hand, it is clear that elements in $\mathcal{RS} (m,\uno,d)$ correspond to vector frames for $\H$.

\pausa
2.
Let $\FS$ be a fusion frame for $\cH\cong \C^d$,
with weights  $w_i>0$ and  subspaces $\ene_i\inc \H$
with $\dim \ene_i = k_i\,$ \fori\,.
Its  fusion frame operator is
$$
S_{\ene_w} = \sum\limits_{i\in \, \IN{m}} w_i^2\, P_{\ene_i} \in \glh^+
$$
(see \cite{[CasKuLi],[CasKuLiRo],GK} for detailed expositions of fusion frames).
Let $U_i\in L(\H,\C^{k_i})$  be a coisometry such that $U_i^*U_i=P_{\ene_i}\,$,
\fori.
Therefore,  the system
 $\cV= \{V_i\}_{i\in \IN{m}} \igdef \{w_i \, U_i\}_{i\in \IN{m}}\,$
 satisfies that $S_{\cV} =S_{\ene_w} \in \glh^+$.
Hence $\cV \in \PRO$ is a projective RS associated to $\ene_w\,$.
Observe that $\cV$ has
the same weights as $\ene_w$ and it
also satisfies that each $\ene_i = R(V_i^*)$.

\pausa
Conversely, given $\RSV \in \PRO$, the sequence $\ene_w = \big(\, \|V_i\| , R(V_i^*) \, \big)_{i\in \, \IN{m}}$
is a fusion frame such that $S_\cV = S_{\ene_w}\,$.
Nevertheless the correspondence is not one to one,
since
any system of coisometries $\{U_i\}_{i\in \, \IN{m}}$ with
$(\ker U_i ) \orto =\ene_i$ produces the same fusion frame $\FS$.
This phenomenon is similar to the correspondence of vector frames
with one dimensional fusion frames.  \EOE
\end{exas}

\begin{rem} \rm In what follows we list some properties and notations about RS's\,:
\ben
\item
Given $\RSV \in \RS$ with
 $S_\cV = \sum_{i\in \, \IN{m}} V_i^*V_i\,$,  then
\begin{equation}\label{hecho1}
\sum_{i\in \, \IN{m}}S_\cV\, ^{-1}\, V_i^*V_i=I_\cH  \ ,
\ \ \text{ and } \ \ \sum_{i\in \, \IN{m}} V_i^*V_i\ S_\cV\, ^{-1}=I_ \cH \ .
\end{equation}
Therefore, we obtain the reconstruction formulas
\begin{equation*}\label{ec recons}
x= \sum_{i\in \, \IN{m}} S_\cV^{-1} \, V_i^* \, (V_i \,x)
= \sum_{i\in \, \IN{m}} V_i^* \,V_i(S_\cV^{-1}\, x) \peso{for every} x\in \cH \ .
\end{equation*}

 \item
For every $\RSV\in \RS$,
we define
the system
$$
\cV^\# \igdef \{V_i\,S_\cV ^{-1}\}_{i\in \, \IN{m}} \in \RS \ ,
$$
called  the
{\bf canonical dual} RS associated to $\cV$.
By Eq. \eqref{hecho1}, we see that
$$
T_{\cV^\#}^* \, T_\cV = \sum_{i\in \, \IN{m}}S_\cV\, ^{-1}\, V_i^*V_i=I_\cH
\peso{and}
S_{\cV^\#} = \sum_{i\in \, \IN{m}}S_{\cV}^{-1}\, V_i^*\, V_i \, S_{\cV}^{-1} = S_{\cV}^{-1}\
\ .
$$

Next we generalize the notion of dual RS's\,: \EOE
\een
\end{rem}

\begin{fed}\label{DV} \rm
Let $\RSV$ and $\RSW \in \RS$.
\ben
\item We   say that $\cW$ is a {\bf dual} RS
for $\cV$ if $ T_\cW^*\, T_\cV = I_\cH\, $,
or equivalently if
$$
\barr{rl}
x & =\suml_{i\in \, \IN{m}} W_i^* \,V_i \,x  \peso{for every} x\in \H  \ .
\earr
$$
\item We denote by  $D(\cV) \igdef \{\cW\in \RS: T_\cW^*\, T_\cV = I_\cH\,\} $,
the set of all dual RS's for a fixed $\cV \in \RS$.
Observe that  $D(\cV) \neq \vacio$ since $\cV^\#\in D(\cV) $.

\EOE
\een
\end{fed}

\begin{rem} \label{MP sost}\rm 
Let $\cV\in \RS$. Then $\cW\in D(\cV)$ \sii  its synthesis
operator $T_\cW^*$ is a pseudo-inverse of $T_\cV\,$.
Indeed,  $\cW\in D(\cV) \iff T_\cW^*\, T_\cV = I_\cH\,$. Observe
that the map $\RS \ni \cW \mapsto T_\cW^*$ is one to one.
Thus,   in the context of RS's   each $(m,\k,d)$-RS has many duals that are $(m,\k,d)$-RS's.
This is one of the advantages of the RS's setting.

\pausa
Moreover,
the synthesis operator  $T_{\cV^\#}^*$ of the canonical
dual $\cV^\#$
corresponds to the Moore-Penrose pseudo-inverse of $T_\cV\,$. Indeed, notice that
 $T_{\cV}\,T_{\cV^\#}^* = T_{\cV}\, S_{\cV}\inv T_{\cV}^*
 \in L(\cK)^+$, so that it
 is an orthogonal projection.
From this point of view, $\cV^\#$ has some optimal properties that come from the theory of pseudo-inverses.  \EOE

\end{rem}

\newcommand{\combi} {\binom}

\section{Optimal erasures and errors}\label{sec erasures}

In coding theory, and using  our terminology,  a signal $x\in \C^d$ is  transmitted encoded via the RS $\RSV$ in the form of $m$ packets $\{V_ix\}_{i\in \IN{m}}$. Then, the loss (or delay) of a number of packets in transmission is treated as if the corresponding components in the RS were ``erased''. Assuming that a number of erasures occurred and taking into account the redundancy of the RS, the reconstruction of the signal could be attempted using a dual RS of $\cV$. The accuracy of this process can be quantified by the norm of an error operator. For example, it is shown in \cite{HolPau} that uniform Parseval frames  (together with their canonical duals) are optimal for 1 erasures in the sense that they minimize a norm of the error operator.
 The authors also prove that equiangular uniform Parseval frames are optimal for 2 erasures (in case that such frames exist). 
In the context of Parseval fusion frames, Bodmann (\cite{Bod}) shows that uniform Parseval fusion frames (a protocol associated to a uniformly weighted projective resolution of the identity, using his terminology) are optimal for 1 erasure, when the dimensions of the subspaces are equal. In addition, he shows that the uniform Parseval protocols that are optimal for multiple erasures have subspaces satisfying the so-called equi-isoclinic condition.

Similar problems have also been considered in different contexts. For example,
in \cite{GK} the authors describe the structure of the optimal frames for one or multiple erasures, but in a non-deterministic setting. Indeed, in this case  a random vector $x$ is estimated from its fusion frame measurements using the Linear Minimum Mean Square Error. 

As it was mentioned in the introduction, our aim is to set the problem studied in  \cite{HanLo} in the context of reconstruction systems. 
That is, given a fixed projective RS $\cV$, the goal is to characterize optimal alternate duals for $r$ erasures. As in \cite{HanLo}, our results imply the existence of a unique dual RS (which is not necessarily projective) that is optimal for the erasure of 1 packet of coefficients. Notice that in this context, optimal alternate duals for $r$ erasures among those duals which are optimal for $r-1$ erasures (see \cite{Bod,Pau,BodPau,HolPau,HanLo,P}). Therefore, by the nature of our results, the optimal duals for $r$ erasures coincide with optimal duals for 1 erasure.

In order to describe the different measures of the reconstruction error 
when an arbitrary packet of coefficients of the fixed RS is erased, we consider the following notions.  Given $j\in \IN{m}$, let
$$
M_j  \in L(\cK)  \peso{given by} M_j \, \big(\,  (y_i)_{i\in \, \IN{m}}\big) =
\big( \, \uno_j \,(i)  \cdot  y_i\,\big)_{ i\in \, \IN{m}}
\ ,
$$
	where $\uno_ j \, : \, \IM \rightarrow \{0,1\}$
denotes the characteristic function of the set $\{j\}\subset\IM\,$. Similarly, we consider the
{\it packet-lost operator} $L_j  \igdef M_{\IM\setminus\{j\}} = I_\cK -M_j\,$.

\pausa
Given $\RSV \in \RS$,
we shall consider a ``blind reconstruction'' strategy in
 case that some coefficient is lost. That is, assuming that the encoded
information $T_\cV \, x\in \cK$ (for some $x\in \mathcal H$) is altered according to the
packet-lost operator $L_j\,$, our reconstructed vector will be
$\hat{x}=T_\cW^*\, L_j \, T_\cV \, (x)$, where
$\RSW\in D(\cV)$
is some dual RS for $\mathcal V$. Then  the reconstruction error will 
be  $x-\hat{x}=x-T_\cW^*\, L_j \, T_\cV \, (x)=T_\cW^*\, M_j \, T_\cV \, (x)
= W_j^*V_j \, x$.

\pausa
In this case, we will use the Frobenius norm 
$\|\cdot\|_{_2}$ to perform the measure of
the operator  $W_j^*V_j \,$. It is worth to note that, unlike the 
vector frame case, $W_j^*V_j $ is not necessary a rank-one operator, 
so its operator (spectral) norm does not coincide with its Frobenius norm. 
We consider that this is a suitable norm to perform the generalization 
of the results in \cite{HanLo} to the RS setting, besides the 
fact that Frobenius norm has nice geometrical properties.
%
%
Consider  the
$m$-tuple
$$
E_1( \cV,\; \cW)=\left( \|I-T_\cW^*\,L_j\,T_\cV \|_{_2}\right)_{j\in \IN{m}}=\left( \|T_\cW^*\,M_j\,T_\cV \|_{_2}\right)_{j\in \IN{m}}=\left(\|W_j^*V_j\|_{_2}\right)_{j\in \IN{m}}
\in \R^m
\ ,
$$
Notice that we can bound uniformly the reconstruction error in terms of the entries of this vector 
for the erasure
of $1$ packet of coefficients (for all $m$ possible choices). In what
follows we shall consider two different measures of the
reconstruction error based on $E_1( \cV,\; \cW)$,
namely the 2-error and the (normalized) worst-case error.

\subsection{Minimizing the 2-error}\label{sec31}
Let $\RSV\in \RS$ 
and let us  denote by
\begin{equation}\label{desc1}
e_1 ^{(2)}(\mathcal V) = \inf_{\mathcal W\in  D(\mathcal V)} 
\|E_1( \cV,\; \cW)\|_{_2}
= \inf_{\mathcal W\in  D(\mathcal V)} 
\left(\sum_{i\in \IM} \|W^*_i \, V_i\|_{_2}^2\right)^{1/2}\ .
\end{equation}
We are interested in the characterization of 
those $\cW\in D(\mathcal V)$ such that 
$\|E_1( \cV,\; \cW)\|_{_2}=e_1^{(2)}(\cV)$.
In other words, we define the set of {\it 1-loss
optimal dual} RS's for $\mathcal V$ with respect to 
$\|E_1( \cV,\; \cdot)\|_{_2}$ as
\begin{equation*}
 D_1^{(2)}(\mathcal V)  \igdef
 \{\cW \in  D(\mathcal V): \|E_1( \cV,\; \cW)\|_{_2}
 =e_1^{(2)}(\cV)\} \ .
\end{equation*}
One of the advantages of this measure of the error is that for 
a projective $\cV$, without any further assumption,
there is always only one dual RS in $D_1^{(2)}(\cV)$, and that this 
system  can be explicitely computed (see Theorem \ref{unico}
and Remark \ref{NO BORR}). 
In particular, using the hierarchies previously 
described, this dual will also be optimal for $r$ erasures.

\begin{rem}\label{no es la MP}\rm
Fix $\RSV \in\PRO$ and $\RSW\in  D(\cV)$. We give a 
matricial interpretation of the 2-error $\|E_1( \cV,\; \cW)\|_{_2} \,$, for the sake of clarity. 
Recall that $\cV \coma \cW \in 
\SP = \bigoplus_{i\in \IM} L(\cH \coma \cK_i)$. Hence they can be 
viewed as block $m\times 1$ column vectors. 
For example, this allows to obtain the equality 
%
$I = T_\cW^*\,T_\cV = \sum_{i\in \IM} W^*_i \, V_i \,$; unfortunately this identity does not allow us to compute the value of $\|E_1( \cV,\; \cW)\|_{_2}$.

\pausa
On the other hand, the 
product 
$T_\cV \, T_\cW^* = \sum_{i,\, j \in \IM} V_i \, W_j^*\in \opk$ is an oblique projector, and 
$\|T_\cV \, T_\cW^*\|_{_2}^2 
= \sum_{i,\, j \in \IM} \|V_i \, W_j^*\|_{_2}^2
= \sum_{i,\, j \in \IM} \| W_j^*\, V_i \|_{_2}^2\,$, since the blocks
 $V_i \, W_j^* \in L (\cK_j\coma \cK_i)$ are the different
``entries" of the $m\times m$ block matrix $T_\cV \, T_\cW^*\,$. 

\pausa
The norm $\|T_\cV \, T_\cW^*\|_{_2}$ is minimal for 
$\cW= \cV^\#$ (since $T_\cV \, T_{\cV^\#}^*$ is selfadjoint). 
Nevertheless, the square of the 2-error  $\|E_1( \cV,\; \cW)\|_{_2}^2 
= \sum_{i\in \IM} \|W^*_i \, V_i\|_{_2}^2 \,$ is 
the square of the norm of the
pinching matrix (i.e., block diagonal truncation) of 
$T_\cV \, T_\cW^*\,$. As we shall see below, its minimum among 
all $\cW\in  D(\cV)$ is not always attained at the Moore Penrose 
pseudoinverse $ T_{\cV^\#}^*\,$ (see Remark \ref{10}). 
\EOE
\end{rem}

\begin{teo}\label{unico}
Let $\RSV\in \PRO$
with weights $\v =(v_i)_{i \in \IN{m}}$. Then
 $D_1^{(2)}(\mathcal V) = \{ \cW_0\}$ i.e.,
there is a {\bf unique} 1-loss optimal dual RS $\cW_0$ for the 2-error. Moreover, if
\ben
\item $D\in L(\cK) $ is the block diagonal matrix $D=\bigoplus\limits _{i\in \, \IN{m}} v_i^{-2}\,I_{\cK_i}\,$, and
\item $S_{\cV, \,D} = T_\cV^*\, D \, T_\cV  = \suml_{i\in \, \IN{m}}  P_{R(V_i^*)} \in \glh^+$
$($\,since $S_{\cV, \,D} \ge (\min\limits_{i\in \IM} \, v_i^{-2}) \cdot S_{\cV}>0$\,$)$,
\een
then the optimal system is $\cW_0 = \{v_i^{-2} V_i \,S_{\cV, D}^{-1}\}_{i\in \, \IN{m}}
 \,$. In particular,
$T_{\cW_0} = D\,  T_\cV \, S_{\cV, \,D}^{-1}\ $.
\end{teo}

\proof
First we check that $\cW_0\in D(\cV)$. Indeed, $T_{\cW_0}^* \, T_\cV = S_{\cV, \,D}^{-1}\,   T_\cV^*\,   D\,  T_\cV  = S_{\cV, \,D}^{-1}\, S_{\cV, \,D} = I$.
Denote by $B_i = v_i^{-2} V_i \,S_{\cV, \,D}^{-1}\,$,
the $i$-th entry of $\cW_0\,$, for every $i \in \IM\,$.
Consider a
dual system
$\RSW \in D(\mathcal V)$.
Since each $V_iV_i^* = v_i^2I_{\cK_i}\,$, then
\beq\label{sale w}
\| W_i^* V_i\|_{_2}^2=\tr(V_i^*W_iW_i^*V_i )
=\tr(W_iW_i^*V_iV_i^*)=v_i^2\,\tr(W_iW_i^*)=v_i^2\,\|W_i^*\|_{_2}^2 \ .
\eeq
In particular, we can compute the 2-error for  $\cW_0\,$:
$$
\|E_1(\mathcal V\coma  \mathcal W_0 )\| ^2 =
\suml_{i\in \, \IN{m}} \|B_i^*\, V_i\|_{_2}^2 \
\stackrel{\eqref{sale w}}{=} \ \suml_{i\in \, \IN{m}} v_i ^2 \|B_i^*\|_{_2}^2
= \suml_{i\in \, \IN{m}} v_i^{-2} \ \| S_{\cV, \,D}^{-1}\ V_i^* \|_{_2}^2
\ .
$$
On the other hand, for every $i \in \IM\,$ we have
\begin{eqnarray}\label{cot 2}
v_i^2\ \|W_i^*\|_{_2}^2&=& v_i^2 \|B_i^* +
(W_i^*- B_i^*)\|_{_2} ^2 \nonumber \\
&=& v_i^{-2}\|S_{\cV, \,D}^{-1}\ V_i^* \|_{_2}^2+v_i^2\ \|W_i^*-B_i^* \|_{_2}^2
+2  \Preal \, \Big(  \,v_i^2\, \tr \,\big[\,   (W_i^*  -B_i^*) B_i\big] \,\Big) \ .
\end{eqnarray}
Let  $t_i= \tr \,\big[\, (W_i^*  -B_i^*) B_i\big] =
v_i^{-2} \tr \big[\, (W_i^*- B_i^* )\, V_i\, S_{\cV, \,D}^{-1}\big] $.
Then we have that
$$
\sum_{i\in \, \IN{m}} v_i^2\ t_i
=  \tr \Big( \ \sum_{i\in \, \IN{m}}(W_i^*- B_i^* )\, V_i\, S_{\cV, \,D}^{-1} \ \Big)
=\tr \big[ \, (T_\cW^*-T_{\cW_0}^*)\ T_\cV \ S_{\cV, \,D}^{-1} \big]=0 \ .
$$
since both $\cW$ and $\cV^\#$ are dual RS's for $\cV$.
Therefore,
summing over $\IM\,$, the third summand of \eqref{cot 2} vanishes  and

\begin{equation}\label{cot 1}
\barr{rl}
\|E_1(\mathcal V\coma  \mathcal W)\|^2  = \suml_{i\in \, \IN{m}}\|W_i^*V_i\|_{_2}^2
&  \ \stackrel{\eqref{sale w}}{=}
 \suml_{i\in \, \IN{m}} v_i^2\ \|W_i^*\|_{_2}^2
\\&\\
&   \ \stackrel{\eqref{cot 2}}{=}
\suml\subim v_i^{-2}\|S_{\cV, \,D}^{-1}\ V_i^* \|_{_2}^2+v_i^2\ \|W_i^*-B_i^* \|_{_2}^2\\&\\
& \geq
\suml_{i\in \, \IN{m}} v_i ^2 \|B_i^*\|_{_2}^2
\stackrel{\eqref{sale w}}{=}
\|E_1(\mathcal V\coma  \mathcal W_0 )\|^2
\ .
\earr
\end{equation}
Therefore $\cW_0 \in D_1^{(2)}(\mathcal V)$. Moreover, if we take another
 $\RSW\in D_1^{(2)}(\mathcal V)$,
then
Eq. \eqref{cot 1} implies that
$\|W_i^*- B_i^* \|_{_2}=0$ for every $i\in \, \IN{m} \, $,
so that   $\mathcal W=\mathcal W_0\,$.\QED

\pausa
We say that a system $\RSV\in \PRO$ is an {\bf uniform} projective RS if
the  weights $\v= (v_i)_{i \in \IN{m}}$ of $\cV$ satisfy that $\v= v \, \uno$ for
some $v>0$.

\begin{rem}\label{NO BORR}\rm
The unique 1-loss optimal dual RS of Theorem \ref{unico}, denoted by $\cW_0\in D_1^{(2)}(\cV)$, can be described in the following way:
If $\v = (v_i)_{i \in \IN{m}}$ are the weights of
the system $\cV\in \PRO$,
consider
$\cU = \{ U_i\}_{i \in \IM} =\{ \frac{V_i}{v_i} \}_{i \in \IM} \,$,
which is a uniform projective RS.
Then
\beq\label{no borr}
S_{\cV, \,D}
= \suml_{i\in \, \IN{m}}  P_{R(U_i^*)}= S_\cU \quad , \quad
\cV = \v \cdot \cU \igdef \{v_i \, U_i\}_{i \in \IM}
\peso{and}
\cW_0 = \v\inv \cdot \cU^\# \ ,
\eeq
because
$\cW_0 = \{v_i^{-2} V_i \,S_{\cV, \,D}^{-1}\}_{i\in \, \IN{m}}
= \{v_i \inv \, U_i \, S_\cU\inv \}\subim \,$.
\EOE
\end{rem}

\begin{cor}\label{cor unif mse}
Let $\RSV\in \PRO$.  Assume that $\cV $ is uniform. 
Then  the unique 1-loss optimal dual RS
for the 2-error  is the canonical dual $\cV^\#$.
\end{cor}
\proof
If the weights of $\cV$ are $\v = v\, \uno$ then, with the notations of
Remark \ref{NO BORR}, we have that $\cV = v \, \cU $,
$S_\cV = v^2 \, S_\cU$  and $\cV^{\#}
= v\inv \, \cU^{\#}$. Then we apply Eq. \eqref{no borr}. \QED

\begin{rem}\label{10} \rm
For most systems $\RSV\in \PRO$ which are not uniform, the unique 1-loss optimal dual RS
for the 2-error $\cW_0\in D(\cV)$ obtained in Theorem \ref{unico} does not coincide with the canonical dual $\cV^\#$. 
For example, it is easy to see that $\cW_0\neq \cV^\#$
if all the weights are different
and if the subspaces $R(V_i^*)$ are not mutually orthogonal. In particular, $\cW_0\neq \cV^\#$ whenever the weights are different and $\cV$ has non zero redundancy. \EOE
\end{rem} 

\subsection{Minimizing the worst-case reconstruction error}\label{wre}

Let $\RSV \in \PRO$.
For each $\RSW \in D(\cV)\,$,
we introduce the worst-case reconstruction error when one packet is lost with
respect to the Frobenius norm. Now, we measure the error vector $E_1(\cV,\; \cW)$ with the maximum of its entries, instead of the euclidean norm used in the previous section. 

\begin{equation}\label{desc}
e_1 (\mathcal V) = \inf_{\mathcal W\in  D(\mathcal V)} \|E_1(\cV,\cW)\|_{_\infty}
= \inf_{\mathcal W\in  D(\mathcal V)} \ \max_{i\in \, \IN{m} \, }  \|T_\cW^* \, M_i \, T_\cV\|_{_2}
= \inf_{\mathcal W\in  D(\mathcal V)} \ \max_{i\in \, \IN{m} \, }  \|W^*_i \, V_i\|_{_2}\ .
\end{equation}
We define the set of {\it 1-loss
optimal dual} RS's for $\mathcal V$ with respect to $\|E_1(\cV,\cdot)\|_{_\infty}$ as
\begin{equation*}
 D_1(\mathcal V)  \igdef\{\cW \in  D(\mathcal V): \|E_1(\cV\coma \mathcal W)\|_{_\infty}=e_1(\cV)\} \ .
\end{equation*}
The study of $D_1(\cV)$ has been considered by Han and L\'opez in \cite{HanLo}
in the particular case of $(m,\uno,d)$-RS's for $\H$, i.e. usual vector
frames. Indeed, since in such case the operators $i\in \IN{m}$, $T_{\cW}^*M_i T_{\cV}$ are rank one operators, their Frobenius and spectral norms coincide, so the WCE coincide with the measure of the error used in \cite{HanLo}.  The use of the Frobenius norm in the definition of the WCE allows to extend naturally the results in \cite{HanLo} to the RS setting using similar techniques.

\begin{rem}\rm 
Recall that  $\RS \inc  \SP= \bigoplus_{i\in \IM} L(\cH \coma \cK_i) \cong \lhk$.
If we fix an  {\bf injective}  $\RSV\in \RSI$, then the map $\|\cdot\|_\cV : \SP \to \R_+$ given by
$$
\|\cW\|_\cV \igdef \|E_1(\cV\coma \cW)\|_{_\infty}  = \max_{i\in \, \IN{m} \, }  \|W^*_i \, V_i\|_{_2}  \peso{for} \RSW\in \SP
$$
is a norm in $\SP$. Indeed, the only non trivial condition is the faithfulness. But
the fact that $\cV\in \RSI$ (i.e.  $V_i$ is surjective for every $i\in\IM$)
assures that $\|\cW\|_\cV= 0 \implies
\|W^*_i \, V_i\|_{_2} = 0 $ \fori $\implies \cW =0$.

\pausa
Since $D(\cV)$ is closed
in $\SP$ with the usual norm and all norms are equivalent on $\SP$, then
 $D(\cV)$ is
$\|\cdot\|_\cV\,$-closed in $\SP$.
Therefore there exist elements $\cW\in D(\cV)$ such that
$\|\cW\|_\cV = \min\limits_{\cN\in D(\cV)} \,
\|\cN\|_\cV  = e_1(\cV)$. Indeed, the intersection of
$D(\cV)$ with a fixed closed ball is a compact set.
On the other hand,
$D(\cV)$ is convex (actually it is an affine
manifold). Since every norm is a convex map,
we have proved the following result: \EOE
\end{rem}

\begin{pro} \label{hay min}\rm
Let $\cV\in \RSI$ be an {\bf injective} system. Then the set
$D_1(\cV)$ of $1$-loss optimal dual RS's for $\cV$ is
non-empty, compact and convex.
\QED
\end{pro}

\begin{teo}\label{GenHL}
Let $\RSV \in \PRO$ with weights $\v =(v_i)_{i \in \IN{m}}$.
If
$$
\| S_\cV^{-1} V_i^*V_i\|_{_2}= v_i^2 \| S_\cV^{-1} P_{R(V_i^*)}\|_{_2} = c  \peso{for every} i\in \, \IN{m} \  ,
$$
then $\cV^\#$, the canonical dual RS of $\cV$, is the unique 1-loss optimal dual RS for $\cV$
(and hence the $r$-loss optimal dual RS for every $r$). In other words, $D_1(\cV) = \{\cV^\#\}$.
\end{teo}

\proof
By Proposition \ref{hay min}, there exists some  $\RSW \in D_1(\cV)$.
Then
$$
\|\cW\|_\cV =  \max_{i\in \, \IN{m} \, } \| W_i^* \, V_i\|_{_2}
\leq
\max_{i\in \, \IN{m} \, } \|S_\cV^{-1}\, V_i^*\, V_i\|_{_2}  = \|\cV^\#\|_\cV =  c \ .
$$
If we denote each $V_iS_\cV^{-1} = C_i\,$, then
$ \| W_i^* V_i\|_{_2}^2 \leq c=  \| C_i^* V_i\|_{_2}^2$ for every $ i\in \, \IN{m} \,$.
Recall that, by Eq. \eqref{sale w},  $\| W_i^* V_i\|_{_2}^2=v_i^2\|W_i^*\|_2^2 \,$,
 since
$V_i V_i^* = v_i^2 \, I_{\cK_i} \,$.
Similarly, we get that each $\|C_i^*V_i\|_{_2}^2=v_i^2\,\|C_i^*\|_{_2}^2$.
Therefore
$
\|W_i\|_{_2}^2\leq \|C_i^*\|_2^2 $ \fori.
Note that
$$
\barr{rl}
\|W_i^*\|_{_2}^2&= \|C_i^*+ (W_i^*-C_i^*)\|_{_2} ^2 = 
\|C_i^*\|_{_2}^2+\|W_i^*-C_i^*\|_{_2}^2+
2 \, \Preal \Big( \, \tr\, \big[ \, (W_i^*-C_i^*)C_i\big] \,\Big)
\earr
$$
and hence $\|W_i^*-C_i^*\|_{_2}^2+ 2 \,
\Preal \Big( \, \tr\, \big[ \, (W_i^*-C_i^*)C_i\big] \,\Big)\leq 0$, for every
$i\in \, \IN{m} \, $. Finally,
$$
\sum_{i\in \, \IN{m}} \tr \, \big[ \, (W_i^*-C_i^*)C_i\big]= \tr \, \Big[ \, (T_\cW^*-T_{\cV^\#}^*)T_\cV S_\cV^{-1} \Big]
=0 \ ,
$$
since both $\cW$ and $\cV^\#$ are dual RS's for $\cV$. Then
$$
0\le  \sum_{i\in \, \IN{m}}  \|W_i^*-C_i^*\|_{_2}^2  =
\sum_{i\in \, \IN{m}} \|W_i^*-C_i^*\|_{_2}^2 + \sum_{i\in \, \IN{m}}
2 \, \Preal \Big( \, \tr\, \big[ \, (W_i^*-C_i^*)C_i\big] \,\Big)  \le 0 \ ,
$$
which implies that $\cW = \{W_i\}_{i\in \, \IN{m}}=\{C_i\}_{i\in \, \IN{m}} = \cV^\#$.\QED

\pausa
A system $\cV \in \RS $ is called a {\bf protocol} for $\mathcal H $  if $S_\cV = I_\cH\,$.
This notion appears in
\cite{Bod}, \cite{P}
(see also \cite{Pau}, where protocols are
related to $C^*$-encodings with noiseless subsystems).

\begin{cor}\label{corpi}
Let $\RSV \in \PRO$
be a projective protocol  for $\H$ (i.e. $S_\cV=I$)
such that  $\|V_i^*V_i\|_{_2} = v_i^2 \, k_i\rai = c$ \fori. Then $D_1(\cV) = \{\cV^\#\}=  \{\cV\}$.
\end{cor}

\proof
By hypothesis $S_\cV =
I_\cH\,$, and hence $\|S_\cV^{-1} V_i^*V_i\|_{_2}
=\|V_i^*V_i\|_{_2}= c$ for every $i\in \, \IN{m} \, $. Thus, the previous
theorem can be applied in this case.

\begin{rem} \rm
Examples of projective protocols as in Corollary \ref{corpi} are
the equi-dimensional uniform projective protocols i.e.,
$\{V_i\}_{i\in \IM}\in \mathcal P(m,k\,\uno,d)$ that are uniform.
These are the analogues of the so-called uniform fusion frames.

\pausa
The $(m,\uno,d)$ case of Theorem \ref{GenHL} is a rephrasing of  \cite[Thm 2.6]{HanLo}, since all vector frames are projective as RS's.  Moreover, we can conclude from the examples \cite[Section 3]{HanLo}
that the optimal dual system $\cW\in D_1(\cV)$ may not be the canonical
dual RS and  may be  not unique for a general $\cV\in \PRO$.
\EOE
\end{rem}

\section{Stability of RS's under erasures of coefficient packets}
\label{rem erasu lb}

In this section we consider a different approach to the erasures problem. 
Indeed given a fixed RS $\cV$, assume that we can identify the set $J\subset \IM$ such that 
$\{V_i\,x\}_{i\in J}$ are the missing or delayed
packets of coefficients in the transmission of the signal $x\in \hil$.
In this case we shall state conditions which assure that the system $\cV_J=\{V_i\}_{i\in \IM\setminus J}$ corresponding to the remaining set of indexes is still a RS. 
Notice that this last fact allows to have perfect reconstruction with suitable dual RS's. Hence, an explicit computation of the canonical dual $(\cV_J)^\#$ is given in this case. 

\pausa
The following statement is a  generalization and a slight  improvement 
of similar results of P. Casazza and G. Kutyniok  \cite{[CasKu2]} 
and Asgari \cite{As} for fusion frames.

\begin{pro}\label{cotas nosiglias }
Let $\RSV \in \RS$  with bounds $A_\cV \coma B_\cV$. 
Fix a subset 
$J\subset \IM$ and consider the matrix 
$M_J \igdef I_d-\sum_{i\in J} V_i^*V_i\,S_\cV^{-1} \in \mat$.   
Then,  
\beq\label{da RS }
\cV_J=(V_i)_{i\in \IM \setminus J}  \quad \mbox{ \rm is a RS for  \ \ $\hil \cong \C^d$ } 
\iff  M_J \in \matinvd  \ .
\eeq
In this case we can compute the following data for $\cV_J\,$: 
\ben
\item The frame operator $S_{\cV_J} = M_J\, S_\cV\,$.
\item The bounds of $\cV_J$ can be estimated by 
$  \frac{A_\cV}{\|M_J^{-1}\|} \le A_{\cV_J}$ and $B_{\cV_J} \le B_\cV\,$.
\item The canonical dual can be characterized as 
$$
(\cV_J)^\#= \{V_i \, S_{\cV_J} \inv \} _{i\notin J} =   
\{\cV^\#_i \, M_J\inv \} _{i\notin J}\igdef
(\cV^\# )_J  \cdot M_J\inv %
   \ .
 $$
 That is, 
 $(\cV_J)^\# $ is the truncation  of  the canonical dual $\cV^\#$ modified with 
 $M_J\inv\,$.
\een
\end{pro}

\proof
It is straightforward to check that $M_J = S_{\cV_J}\,S_\cV^{-1} \implies 
S_{\cV_J} = M_J\, S_\cV\,$. 
This last fact
implies the equivalence of Eq. \eqref{da RS }. 
On the other hand, 
$$A_\cV\, \|M_J^{-1}\|^{-1} \leq \|(M_J\, S_\cV)^{-1}\|^{-1}= \|S_{\cV_J}^{-1}\|^{-1}  
= A_{\cV_J} \ .
$$ 
The fact that $0<S_{\cV_J} \le S_\cV$ assures that  
$B_{\cV_J} \le B_\cV\,$.
Let us  denote 
$\cV^\#  = \{W_i \}\subim $ and 
$(\cV^\# )_J = \{W_i \} _{i\notin J}\,$. Then the formula 
$S_{\cV_J} = M_J\, S_\cV$  
gives the equality  
\beq
(\cV^\# ) _J  \cdot M_J\inv \igdef 
 \{W_i \, M_J\inv \} _{i\notin J}= 
 \{V_i\, S_\cV\inv  \, M_J\inv \} _{i\notin J}
 = \{V_i \, S_{\cV_J} \inv \} _{i\notin J} = (\cV_J)^\#   \ .
 \QEDP
\eeq

\begin{rem}\rm
In \cite{[CasKu2]} P. Casazza and G. Kutyniok stated the 
sufficient (and easily computable) condition 
$\suml_{i\in J} \|V_i\|_{sp}^2 
<A_\cV\,$ for the invertibility 
of the matrix $M_J$ of Proposition \ref{cotas nosiglias }. 
Indeed,  if $\|\suml_{i\in J} V_i^*\, V_i\|_{sp} < A_\cV\,$   
(compare with the condition $\suml_{i\in J} \|V_i\|_{sp}^2<A_\cV\,$),  
then 
$\|I_d-M_J\|_{sp}<1$ and $ M_J \in \matinvd$. 
They also give the  estimation 
$ A_\cV-\suml_{i\in J} \|V_i\|_{sp}^2\le A_{\cV_J}\,$. This follows from 
$$
\barr{rl}
A_\cV-\sum_{i\in J} \|V_i\|_{sp}^2 & \leq  \ 
A_\cV-\|\sum_{i\in J} V_i^*\, V_i\|_{sp} \le 
\frac{A_\cV}{\|M_J^{-1}\|_{sp}} 
\stackrel{Prop. \ref{cotas nosiglias }}{\le} A_{\cV_J}\ .
\earr
$$ 
The result of Asgari in \cite{As} is similar to Eq. \eqref{da RS }, 
but stated for fusion frames and assuming that $|J|=1$, with a 
different lower bound for $\cV_J\,$. \EOE
\end{rem}

\section{Approximation by projective RS's}\label{Approx}
\pausa
Given a fixed $\cV \in \PRO$, notice that the optimal dual RS's obtained in Theorems \ref{unico}, and \ref{GenHL} ,  are not projective RS's, in general. Although there could be some projective elements in $D(\cV)$
(we shall focus this problem in the following section),
we are interested in
those $(m,\k,d)$-projective RS's
that are closest, with respect to
some distance function, to a fixed $\RSS \in D(\cV )$ which has some desired
properties. Given $\cW\in \RS$,
we consider
$$ d(\cS \coma \cW)\igdef \|T_{\cS}-T_{\cW}\|_{_2}= \|T_{\cS}^*-T_{\cW}^*\|_{_2} \ ,
$$
the distance between their synthesis (or analysis) operators.
Hence, we seek for  $\cW_0 \in \PRO$ that minimize
$d(\cS \coma \cW)$ among the projective RS's.
In what  follows we will describe the structure of such (unique) minimizers in case $\cS$ is an injective RS.
As one would expect, its ``directions" are the coisometries of the
polar decompositions of the coordinate operators $S_i $ of $\cS$,
while its weights are the ``averages" of their singular values.
We need first some preliminary results:

\pausa
Given $k, n \in \N$ such that $k\le n$, we denote by
$$\mathcal I(k\coma n) = \{U \in L(\C^k \coma \C^n) : U^*U = I_k\} \ ,
$$
the set of  isometries from $\C^k$ into $\C^n$. The following result 
can be found in \cite{Bha}: 
\begin{lem}\label{lem aprox} \rm
Let $k, n \in \N$ such that $k\le n$, and  let $A\in \cM_{n,k}(\C)$ be a full rank matrix with
polar decomposition $A=U\,|A|$ with $U\in \mathcal I(k\coma n)$.
Then
$\|A-U \|_{_2} \ = \ \min\limits_{V\in \mathcal I(k\coma n)} \ \|A - V\|_{_2} \ 
$. \QED
\end{lem}
Recall that for every $A\in \cM_{n,k}(\C)$ its polar decomposition satisfies
\beq\label{UA AU}
A = U|A| = |A^*| U  \implies U^*A = |A| \peso{and} A^* = U^* |A^*|  \ ,
\eeq
where $U \in \cM_{n,k}(\C)$ has $\ker U= \ker A$.
As expected, the aproximation of a RS by a projective RS $\{S_i\}_{i\in \IN{m}}$ relies on the coisometries that best approximate each $S_i$, which are determined by the polar decomposition of $S_i$. Nevertheless, one has to determine the appropiate weights which combine those coisometries. That is done in the following Proposition.

\begin{pro}\label{distancia desde los ff}
Let $\RSS\in \RSI$.
Then there exists a unique
\begin{equation}\label{ecua min}
\cW_0 \in \PRO \peso{\rm such that}
d(\cS \coma \cW_0)=\min _{\mathcal W\in \PRO}  d( \cW\coma \cS)\ ,
\end{equation}
and it is given by $\mathcal W_0=\{\alpha_i\ U_i\}_{i\in \, \IN{m}}\ $
where each 
$\displaystyle \alpha_i
=\frac{ \tr |S_i|}{k_i} \, $ and 
 $S_i =U_i \,|S_i |$ is the polar decomposition of each $S_i\,$.
 \end{pro}
\proof
Let $\cW\in\PRO$ be a system such that the minimum in \eqref{ecua min} is attained
at $\cW$.
Denote by $\w =  (w_i)_{i \in \IN{m}}\in \R_+^m $ the weights of $\cW$. Notice that
$$
\|T_{\cS}^*-T_\cW^*\|_{_2}^2
=\sum_{i=1}^m \|S_i^*-W_i^*\|_{_2}^2
\peso{and \ \ each} W_iW_i^*=w_i^2\,I_{\cK_i} \ ,
$$
Thus, each isometry $w_i^{-1}\,W_i^*\in \cI(k_i\coma d)$ attains the minimum in the optimization problem
$$
\big\| \, w_i^{-1}\,W_i^* - \frac{S_i^*}{w_i} \, \big\|_{_2}
=
\min_{X \in \cI(k_i\coma d)} \ \big\| \, X-
\frac{S_i^*}{w_i} \, \big\|_{_2}
 \ ,
$$where, by hypothesis, each $w_i^{-1}\,S_i^*$ is a full rank linear transformation.
By Lemma \ref{lem aprox} we get that $w_i^{-1}\,W_i^*=U_i^*\, $,
the isometry of the polar decomposition
$$
w_i^{-1}  S_i ^* = U_i^* \ |w_i^{-1} \,  S_i ^* |
= w_i^{-1} \ \big(\, U_i^* \ |S_i ^* | \, \big)\
\stackrel{\eqref{UA AU}}{\implies}  \ S_i = U_i \  |S_i |  \peso{and}
\ker U_i = \ker S_i
\ ,
$$
and hence $W_i=w_i\,U_i\,$.
Next we  show that each $w_i =  \frac{ \tr |S_i|}{k_i} \ $.
Fix $i\in \IM\,$. Then
\[
\|S_i  - w_i  \,U_i\|_{_2} \
= \ \min_{\alpha>0} \   \|S_i  - \alpha \,U_i\|_{_2} \ .
\]
Therefore $w_i \,\cdot \|U_i\|_{_2}$ is the norm of the orthogonal projection of
$S_i $  to the line $\R\, U_i\,$, using the $\R$-inner product
$\pint{A,B}=\Preal \big[ \tr (B^*A)\big]$.
It can be computed explicitly:
$$
0\le \frac{\tr \, |S_i|}{\|U_i\|_{_2}}
\ \stackrel{\eqref{UA AU}}{=} \
\frac{\tr \, \big(\, U_i^* \, S_i  \big)}{\|U_i\|_{_2}}
= \left|\, \pint { S_i \, , \, \frac{U_i}{\|U_i\|_{_2}}}\,  \right|
= \|P_{ \R\, U_i} ( S_i )\|_{_2} =
w_i  \, \cdot \|U_i\|_{_2}
\ , $$
for every $i \in \IM\,$. 
Then we obtain the equalities  $\displaystyle
w_i = \|U_i\|_{_2}^{-2} \, \tr |S_i|  = \frac{ \tr |S_i|}{k_i} \ $.
\QED

%
%

\def\GV{\cV(G\coma V)}

\section{Examples}\label{Opti}
In this section we present a variety of examples related with the previous sections. First, we exhibit a family of (projective) RS's which satisfies the hypotheses of Corollary \ref{cor unif mse} and Theorem \ref{GenHL}, the so-called group RS. In particular, the canonical dual of a group RS  (which is also a group RS) is the unique 1-loss optimal dual RS for the 2-error as well as for the WCRE.
The remaining subsections are devoted to the study of particular cases of RS where  $D(\cV) \cap \PRO \neq \vacio \ .$  
The first examples show that for certain projective systems we can explicitly construct projective duals, which in general will not coincide with the canonical duals.  
The last example describes a Riesz RS whose unique dual RS (i.e. the canonical dual) is not projective. This leads to a characterization for Riesz RS`s with projective canonical dual.

\subsection {Group reconstruction systems}
We begin by rephrasing some basic notions and results from \cite{HanLar} in the RS's setting:
Let $\cK_0\cong \C^k$ and $V  \in L(\cH \coma \cK_0)$.
Given a unitary representation
$G\ni g \mapsto U_g \in \uni$ of a finite group $G$ in $\cU(\cH)$ we say that
\[\GV \igdef \{V \,U_g\}_{g\in G}\]
is a  $G${\bf-reconstruction system}  (shortly,  $G$-RS) if
$\GV  \in \cR\cS (m, k\uno, d)$, where $m=|G|$. In this case the space $\cK = \cK_0^m\,$. 
If $\ese = V^*(\cK_0)$, this is equivalent to the fact that
$$
\gen \Big\{ \ \bigcup \limits_{g\in G} \ U_g (\ese)\, \Big\} = \cH\ ,
$$
where $V\in L(\cH \coma \cK)$ is {\bf the base operator} for $\GV$.
This definition of G-RS reduces to that of $G$-frame
in the  vector frames setting. Following  \cite{HanLar} and 
\cite{Han} we state a series of properties of
G-RS's whose proofs are similar to the frame case: Fix  $V\in L(\cH \coma \cK_0)$. 
\ben
\item 
Observe that the system $\GV$ is:
\ben
\item Projective (and uniform) if in addition $VV^*=v^2 I_{\cK}$ for some $v>0$ ;
\item Injective if in addition $V^*$ is injective, in which case also $R(V\,U_g) = \cK$ for every $g\in G$.
\een
\item
Notice that
the RS-operator of $\GV $ has the following structure:
\[
S_{G\coma V} \igdef S_{\GV} = \sum\limits_{g\in G} U_g^* \, V^* \, V\,U_g
= \sum\limits_{g\in G} U_{g\inv} \, V^* \, V\,U_g  \ .
\]
\item
The RS-operator $S_{G\coma V}$ (and therefore $S_{G\coma V}^{-1}$)
commutes with the unitary representation of $G$:
\beq\label{kettle}
U_{h } \cdot \, S_{G\coma V} =  S_{G\coma V}\cdot U_{h } 
\peso{for every} h \in G\ .  
\eeq
\item In particular,
the canonical dual of a $G$-RS is another $G$-RS:
\[
\GV ^\# = \{V\,U_g\, S_{G\coma V} ^{-1}\}_{g\in \, G}=
\{V\,S_{G\coma V} ^{-1} \,U_g\,\}_{g\in \, G}  = \cV(G\coma V\,S_{G\coma V} ^{-1} )\ .
\]
\een
In order to apply our previous results, assume now that the base operator $V\in L(\cH,\cK)$ satisfies $VV^*=v^2\, I_\cK$. Then $\GV$
is an equi-dimensional uniform projective RS. Hence, 
Corollary \ref{cor unif mse} implies that the
$1$-loss optimal dual RS for $\GV$ for the 2-error is its canonical dual $\GV^\#$.
In this case 
then $\GV$ also satisfies the hypothesis of Theorem \ref{GenHL}: 
$\GV $ is a $G$-projective RS, and 
\[
\| (S_{G\coma V}^{-1} \,U_g^*\, V^*)(V \, \,U_g)\|_{_2}=
 \| \,U_g^* \, (S_{G\coma V}^{-1} \,V^*V )\,U_g\, \|_{_2}
=  \| S_{G\coma V}^{-1} \,V^*V\|_{_2}= c    \ . 
\]
Thus, the canonical dual
of such $G$-RS's is the unique 1-loss optimal dual for the worst-case error.
If the base operator $V\in L(\cH,\cK)$ is surjective, so that $\GV$ is an injective RS then, using Proposition \ref{distancia desde los ff}, the projective
RS nearest to $\GV ^\# $ can be computed in the following way:
For every $g \in G$, we have that
$$
|V \, S_{G\coma V}^{-1} \,U_g |^2 = U_g^*\,S_{G\coma V}^{-1}\,V^*\,V\,S_{G\coma V}^{-1}\,U_g =
U_g^*\,|VS_{G\coma V}^{-1}|^2 \,U_g
\ .
$$
Taking square roots at both sides, we get that $|V\,U_g \, S_{G\coma V}^{-1}| =
U_g^*|VS_{G\coma V}^{-1}| U_g $ for every $g\in G$.
Therefore, if we consider the polar decomposition
$V \, S_{G\coma V}^{-1}= W|V \, S_{G\coma V}^{-1}|$ of $V \, S_{G\coma V}^{-1}\,$, then also
$$
V \, S_{G\coma V}^{-1} \,U_g
= \big(  \, W \, U_g \big) \ \big( U_g^*\, |VS_{G\coma V}^{-1}| U_g \, \Big)
= \big(  \, W \, U_g \big) \  |V\,  U_g \, S_{G\coma V}^{-1} |
$$
is the polar decomposition of each entry $V \, S_{G\coma V}^{-1} \,U_g \,$ of $\GV^\#$.
In conclusion, if we denote
$$
w=  \frac{\tr \, |V \, S_{G\coma V}^{-1}|}{k}\  \ ,  \peso{then}
\cV(G\coma w \, W) = \big\{ \, w\, W \, U_g \,\}_{g\in G}
$$
is the best projective approximation of $\GV^\#$
provided by Proposition \ref{distancia desde los ff}.
It is clear from the previous computations that 
it is again a $G$-RS.

\subsection{Dual projective systems}
Next we consider an example of a system $\RSV\in \PRO$ with projective dual
systems but such that $\cV^\# \notin \PRO$\,:

\begin{exa}\label{ej por venir}\rm 
Let $d =3$, $m=2$ and $\k = (2,2)$. Let $V_1 $ and $V_2 \in L(\C^3 \coma \C^2)$ be given by
$$
V_1 (x,y,z) = (y,z) \peso{and} V_2 (x,y,z) = (x, z) \peso{for every} (x,y,z) \in \C^3 \ .
$$
Then $\cV = (V_1 \coma V_2) \in \PRO$ with weights $\uno_2\, $. If
$\ese_1  = \{e_1\}\orto $
and $\ese_2  = \{e_2 \}\orto $, then
$$
S_\cV = V_1^* \, V_1 + V_2^* \, V_2 = P_{\ese_1} + P_{\ese_2} = \bm{ccc} 1&0&0\\0&1&0\\0&0&2\em  \ .
$$
Therefore $\cV^\#\notin \PRO$ since 
$S_\cV \inv \, V_1^* (u,v) = (0, u, \frac v2 \, )$ for $(u,v)\in \C^2$,
so that the entry $V_1 \, S_\cV\inv $ of $\cV^\#$ is not a multiple of a coisometry.

\pausa
Let $\cW = (W_1 \coma W_2) \in \PRO$ and assume that  $T_\cW^* \, T_\cV = W_1^* \, V_1 + W_2^* \, V_2 = I_3\,$.
Denote by $\{v_1 \coma v_2\}$ the canonical basis  of $\C^2$ and by $\{e_1 \coma e_2\coma e_3\}$
that  of $\C^3$.
Then, easy computations using the definition of $\cV$ show that
\beq\label{miercoles}
e_3 = W_1^* \, v_2 + W_2^* \, v_2  \quad  , \quad  e_2 = W_1^* V_1 \, e_2 = W_1^* \,  v_1
\peso{and} e_1 = W_1^* V_2 \, e_1 = W_2^* \, v_1 \ .
\eeq
The last two equalities show that both $W_1^*$ and $W_2^*$ should be isometries with weight $1$. But in this case
$\|W_1^* \, v_2\| = \|W_2^* \, v_2\| = 1$ and their sum also has norm one.
Let $\omega = \frac 12 + i \frac {\sqrt 3}{2} \,$. Then $|\omega| = 1$ but $\omega+\ov{\omega} = 1$.
Then we can define $W_1^* \coma W_2^* \in L(\C^2 \coma \C^3)$ by
$$
W_1^*(x, y) = (0,x, \omega\, y)  \peso{and} W_2^*(x, y) = (x,0, \ov{\omega}\, y) \peso{for every} (x,y) \in \C^2 \ .
$$
The $W_i^*$  are isometries and satisfy the three conditions of \eqref{miercoles}. Therefore,  the system
$\cW = (W_1 \coma W_2)$ lies in  $\PRO$ and it   is a dual-RS for $\cV$.

\pausa
Nevertheless, if we  consider $V_1$ and $V_2$ as operators in $L(\R^3 \coma \R^2)$,
 then  such a $\cW$ can not exists in the real case. Indeed,
looking at Eq.  \eqref{miercoles}, we can deduce that
$W_1^* \,v _2  \in \{e_2\}\orto $ and  $W_2^* \,v _2  \in \{e_1\}\orto $.
These facts, together with the
equality  $e_3 = W_1^* \, v_2 + W_2^* \, v_2$ imply that both $W_1^* \,v _2  $
and $W_2^* \,v _2 \in \gen\{e_3\}$,
which is impossible in the real case.
\EOE
\end{exa}

\begin{exa}\rm 
We can generalize Example \ref{ej por venir} in the following way:
Assume that the projective system
$\RSV \in \PRO$ has the property that all the projections $P_i = P_{R(V_i^*)}$ are pairwise commuting.
In this case $
D(\cV) \cap \PRO \neq \vacio \ .
$

\pausa
Indeed, suppose first that all the weights of $\cV$ are $1$.
Then $S_\cV = \sum_{i\in \IM} \, P_i \ $.
The commutation hypothesis assures that, by taking all the possible intersections among the ranges of the projections $P_i\,$, we get a 
family of projections $(Q_j )_{j \in \In}$
such that
\ben
\item $Q_i \, Q_j = 0$ if $i\neq j$ and 
$\sum_{j\in \In} \, Q_j = I_\H\,$.
\item $S_\cV =  \sum_{j\in \In} \, r_j \, Q_j \,$ with $r_j \in \IM$ 
for every $j \in \In\,$.
\item For every $i \in \IM$ there exists $\mJ _i \inc \In $ such that $P_i =  \sum_{j\in \mJ_i} \, Q_j \ $.
\een
We  construct  the system $\RSW \in D(\cV) \cap \PRO $ as follows:
let $W_i = V_i \, U_i\,$, where 
$$
U_i =  \sum_{j\in \mJ_i} \, \eps_{ij} \, Q_j  \peso{for some} 
\eps_{ij} \in   \{\ 1 \ , \ \ -1 \ , \ \  \omega \ , \ \
\ov{\omega} \ \} \ \ \ , \
\peso{where} \omega = \frac12 + i \, \frac{\sqrt 3}{2} \  .
$$
Note that since all $|\eps_{ij}|=1$, then $\cW \in  \PRO $. 
A careful selection of these coefficients, taking account the 
parity of the numbers $r_j = |\,\mS_j\,|$, where $\mS_j = \{i \in \IM : j \in \mJ_i\}$ for  $j \in \In\,$, allows to find such a
$\cW$ such that  $\cW\in D(\cV)$.
%
%
The general case follows from the previous case. Indeed,
if $\cV$  has weights $\v = (v_i)_{i\in \IM}\,$, we replace the previous
$\cW$ by $\cW_\v = \{v_i^{-2} \, W_i\}_{i\in \IM}\,$.  \EOE
\end{exa}

\begin{rem}[Projective dual pairs]\rm The previous example gives a method to construct
pairs of projective RS's $(\cV,\cW)$ such that $\cW\in D(\cV)$.
Moreover, this method shows that for every choice of parameters $(m,\k,d)$ such that
$\sum_{i\in \IM}k_i\geq d$ there exist $\cV,\,\cW\in \PRO$ such that $\cW\in D(\cV)$.
Indeed, to find such a {\bf projective dual pair} $(\cV,\cW)$ we construct
$\cV\in \PRO$ in such a way that the projections
$P_i=V_i^*V_i$ for $i\in \IM$ are pairwise commuting
(i.e. that are simultaneously diagonalizable by an orthonormal basis of $\cH$).
Then, we can apply the construction in the example above
to obtain explicitly the projective dual $\cW$.

\pausa
These facts show that projective dual pairs are indeed more frequent than
projective protocols i.e. RS's $\cV$ such that $(\cV,\cV)$ is a projective dual pair,
since it is known that there are choices of parameters $(m,\k,d)$ for
which no projective $(m,\k,d)$- protocols exist (see \cite[Example 3.1.2.]{MRS}).
\EOE
\end{rem}

\subsection{Riesz reconstruction systems}
\pausa
The elements of $\RS$ are called Riesz RS's if $\dim \cK
= \tr \k
= d = \dim \cH\,$.
In this case, every $\RSV\in \RS$ has the property that both the synthesis operator $T_\cV$
and the analysis operator $T_\cV^*$ are invertible.
Also $T_{\cV^\#}^* =  S_\cV\inv \,T_\cV^* = T_\cV\inv\,$, and
\beq\label{sum dir}
\cH \ = \ R(V_1^*) \oplus \dots \oplus R(V_m^*)  \ \ \ \
\mbox{(direct sum, but not necessarily orthogonal)}\ ,
\eeq
because the sum gives always $\cH$ but in this case
the sum must be direct by dimensional reasons.
On the other hand, if $\RSV\in \RS$ is a Riesz RS for $\cH$ then
\beq\label{RRS}
D(\cV) = \{\cV^\#\}  \ ,
\eeq
since
the  only left inverse of $T_\cV$ is $T_\cV^{-1}=T_{\cV^\#}^*$.
Recall that $T_\cW^*\, T_\cV = I_\cH$ for  every $\cW\in D(\cV)\,$, and
that the map $\RS \ni \cW \mapsto T_\cW^*$ is one to one.

\begin{exa}\label{Directa}\rm 
Let $d=4$, $m=2$ and $\k = (2,2)$. We now construct a (necessarily) Riesz $\cV\in \PRO$ such that $D(\cV)\cap \PRO = \vacio$. Let $V_1\coma V_2 \in L(\C^4 \coma \C^2)$ be given by
\beq\label{V1V2}
\barr{rl}
V_1(x_1\coma x_2\coma x_3\coma x_4) & = (x_1\coma x_2) \peso{and} V_2 (x_1\coma x_2\coma x_3\coma x_4) = (x_3\coma \frac{x_2 -x_4}{\sqrt 2} ) \ , \earr
\eeq
for $(x_1\coma x_2\coma x_3\coma x_4)\in \C^4$. It is easy to see that
$\cV = (V_1 \coma V_2)\in \PRO$ with weights $(1,1)$.
Let us denote by $S = \ker V_2 = \gen\{e_1 \coma e_2+e_4\} \inc \C^4$.
Given $\cW\in D(\cV)$, the equality
$$
 W_1^*\, V_1 + W_2^*\, V_2 = I_\cH
 $$
implies that $W_1^* \, V_1 \, x = x $ for every $x\in S$.  Then $W_1^* \in
L(\C^2 \coma \C^4)$ is completely determined as the inverse of $V_1\big|_S : S \to \C^2$.
But we have that
\bce
$\|V_1 \, e_1\| = \|e_1\| = 1$  \ \ while  \ \
$\|V_1(e_2 +e_4) \| = \|e_2\| = \ds \frac{\|e_2+e_4\|}{\sqrt 2} \ $.
\ece
Then $V_1\big|_S$ is not a multiple of an isometry and neither is $W_1^*\,$.

\pausa
We can enlarge the previous example in order to get a RS with redundancy and without projective duals. Indeed, consider 
$\cV_0=(V_1 \coma V_2 \coma V_3 ) \in 
\cP (3\coma (2,2,2) \coma 4)$ obtained from $\cV$ by adding any 
coisometry $V_3 \in L(\C^4\coma \C^2)$ such that also $\ker V_3 = \ese$. Then, arguing as before, we conclude that there is no  
$\cW = (W_1 \coma W_2 \coma W_3 ) \in D(\cV_0)$ 
such that $W_1$ is a multiple of a coisometry. 
\EOE
\end{exa}

\begin{rem}\rm Assume that $\RSV \in \RS$ is a Riesz RS. Then, arguing as in the previous example, it is easy to see that the following conditions are equivalent:
\ben
\item $D(\cV)\cap \PRO \neq \vacio$.
\item $\cV^\# \in \PRO $.
\item If we denote by $S_i = \bigcap\limits_{j\neq i}  \ker V_i
= \big( \bigoplus \limits_{j\neq i} R(V^*) \, \big) \orto $,  then $V_i \big|_{S_i} \in L(S_i\coma \cK_i)$
is a multiple of an isometry,
for every $i \in \IM\,$.
\een
These conditions are fulfilled if the sum of Eq. \eqref{sum dir}
is orthogonal. Also if every $k_i = 1$. But there exist other cases in which  $D(\cV)\cap \PRO \neq \vacio$. For example, if we
take the  operator $V_1 $  of  Eq. \eqref{V1V2}, and consider
$$
V_3 (x_1\coma x_2\coma x_3\coma x_4) = (x_1- x_3\coma x_2 -x_4) \peso{for}
(x_1\coma x_2\coma x_3\coma x_4) \in \C^4 \ ,
$$
then the condition 3 is satisfied by
$\cV' = (V_1 \coma V_3)\in \cP(2,(2,2), 4)$.
\EOE
\end{rem}

\fontsize {9}{9}\selectfont

\end{document}